%% Authors:     Gordon and Luecke
%% Title:       Toroidal and Boundary-Reducing Dehn Fillings
%% Accepted:    Topology and Its Applications           (Code:  TOPOL1281)
%% VERSION:     10/7/97
%% Initialized: 12/5/96
%% Report-no:	ut-ma/980001
%%%%%%%%%%%%%%%%%%%
\input amstex
\documentstyle{amsppt}
\magnification=\magstep1
\pageheight{9truein}
\pagewidth{6.5truein}
\NoBlackBoxes
\input epsf
%%%%%%%%%%%%%%%%%%%
\def\que{{\Bbb Q}}
\def\real{{\Bbb R}}
\def\zed{{\Bbb Z}}
\def\bq{\text{\bf q}}
\def\nhd{\operatorname{nhd}}
\def\C{{\Cal C}}
\def\D{{\Cal D}}
\def\chix{\raise.5ex\hbox{$\chi$}}

\def\fig#1#2#3{\vbox{\centerline{\epsfysize#1truein \epsfbox{#2}}
		     \centerline{\smc Figure #3} }}
\def\oLambda{\overline{\Lambda}_0} 
%%%%%%%%%%%%%%%%%%%
\topmatter
\title Toroidal and Boundary-Reducing Dehn Fillings\endtitle
\author C. McA. Gordon and J. Luecke\endauthor
\leftheadtext{C. M{\fiverm c}A. GORDON and J. LUECKE}
\address Department of Mathematics, The University of Texas at Austin,
Austin, TX 78712-1082\endaddress
\email gordon\@math.utexas.edu; luecke\@math.utexas.edu\endemail
\abstract  
Let $M$ be a simple 3-manifold with a toral boundary 
component $\partial_0 M$. 
If Dehn filling $M$ along $\partial_0 M$ one way produces a toroidal 
manifold and Dehn filling $M$ along $\partial_0 M$ another way produces a 
boundary-reducible manifold, then we show that the absolute value of the 
intersection number on $\partial_0 M$ of the two filling slopes is at most two. 
In the special case that the boundary-reducing filling is actually a solid 
torus and the intersection number between the filling slopes is two, more 
is said to describe the toroidal filling.
\endabstract  
\endtopmatter

\document 

\head  1. Introduction\endhead 

Following \cite{W2}, let us call a compact, orientable 3-manifold 
{\it simple\/} if it contains no essential sphere, disk, torus or annulus. 
Let $M$ be such a manifold, and let $\partial_0 M$ be a torus component 
of $\partial M$. 
If $\gamma$ is a {\it slope\/} on $\partial_0M$ (the isotopy class of an 
essential unoriented circle), then as usual $M(\gamma)$ will denote the 
manifold obtained by {\it $\gamma$-Dehn filling\/} on $M$. 
Thus $M(\gamma) = M\cup V_\gamma$, where $V_\gamma$ is a solid torus,  
glued to $M$ by a homeomorphism from $\partial_0M$ to $\partial V_\gamma$ 
taking $\gamma$ to the boundary of a meridian disk of $V_\gamma$. 

If $\gamma,\delta$ are two slopes on $\partial_0M$ such that $M(\gamma)$ 
and $M(\delta)$ are not simple, then there are several results giving upper 
bounds for $\Delta (\gamma,\delta)$ (the minimal geometric intersection 
number of $\gamma$ and $\delta$) for the various possible pairs of essential 
surfaces that arise, which in many cases are best possible (see \cite{W2} 
for more details). 
In the present paper we dispose of one of the remaining cases, namely, that 
in which the surfaces in question are a torus and a disk. 

\proclaim{Theorem 1.1} 
Let $M$ be a simple 3-manifold such that $M(\gamma)$ is toroidal and 
$M(\delta)$ is boundary-reducible. 
Then $\Delta (\gamma,\delta)\le2$.
\endproclaim 

Examples showing that this bound is best possible are given in \cite{HM1} 
and \cite{MM}. 

It was previously known that $\Delta (\gamma,\delta)\le 3$, by \cite{W1}. 
We shall therefore assume from now on that $\Delta (\gamma,\delta)=3$, and 
eventually obtain a contradiction. 
We shall also assume that $\partial M$ has exactly two components 
$\partial_0M$ and $\partial_1M$, each a torus, and that $M$ is a 
$\que$-homology cobordism between them, since otherwise 
$\Delta (\gamma,\delta) \le1$ by \cite{W2, Theorem~4.1}. 

In \cite{GLu3} it was shown that if $M(\gamma)$ is toroidal and 
$M(\delta) \cong S^3$ then $\Delta (\gamma,\delta)\le2$, and much of the 
proof of Theorem~1.1 consists of carrying over the arguments there to the 
present context. 
(In \cite{GLu3} we analyzed the intersection of the punctured surfaces in $M$ 
coming from an essential torus in $M(\gamma)$ and a Heegaard sphere in 
$M(\delta)$; here the Heegaard sphere is replaced by an essential disk.) 
This quickly leads to a proof of Theorem~1.1 when $t$, the number of points 
of intersection of the essential torus in $M(\gamma)$ with the core of 
$V_\gamma$, is at least~4. 
This is completed in Section~4. 
Sections~2, 3 and 4 parallel the corresponding Sections of \cite{GLu3}; 
(in Section~2 it is also pointed out that Sections~6 and 7 of \cite{GLu3} 
carry over without change). 
The main divergence in the proofs is in  the case $t=2$. 
In \cite{GLu3} this case was handled by showing that the associated knot 
$K$ in $S^3$ was strongly invertible, and then appealing to a result 
of Eudave-Mu\~noz \cite{E-M}. 
Since no analog of this result is available in our present setting, 
we have to argue this case directly. 
This is done in Section~5. 

In Section 6 we specialize to the case that the boundary-reducible filling 
of $M$ is a solid torus. 
Infinitely many examples of simple 3-manifolds $M$ and 
slopes $\gamma,\delta$ with $\Delta (\gamma,\delta)=2$, $M(\gamma)$ 
toroidal, and $M(\delta)$ a solid torus, are given in \cite{MM}.  
In these  examples, there is a twice-punctured essential torus in $M$ whose 
boundary slope is $\gamma$. 
Section 6 is devoted to showing this always happens: 

\proclaim{Theorem 6.1} 
Let $M$ be a simple 3-manifold such that $M(\gamma)$ is toroidal, 
$M(\delta)$ is a solid torus, and $\Delta (\gamma,\delta)=2$.  
Let $K_\gamma$ be the core of the attached solid torus in $M(\gamma)$ 
and $\hat T$ be an essential torus in $M(\gamma)$ 
that intersects $K_\gamma$ minimally. 
Then $|K_\gamma \cap \hat T| =2$.
\endproclaim 

We thank the referee for his helpful comments, and in particular for 
suggesting a simpler proof of Lemma~5.9. 

\head 2. The graphs of intersection\endhead 

Let $K_\gamma (K_\delta)$ be the core of $V_\gamma$ (resp. $V_\delta$). 
Let $\hat T$ be an essential torus in $M(\gamma)$. 
We assume that $\hat T$ meets $K_\gamma$ transversely, and that $t= 
|\hat T\cap K_\gamma|$ is minimal (over all essential tori $\hat T$ in 
$M(\gamma)$). 
We may also assume that $\hat T\cap V_\gamma$ consists of meridian disks, 
so $T= \hat T\cap M$ is a punctured torus, with $t$ boundary components, 
each having slope $\gamma$ on $\partial_0M$. 

Similarly, let $\hat Q$ be an essential disk in $M(\delta)$, meeting 
$K_\delta$ transversely, with $q= |\hat Q\cap K_\delta|$ minimal over all 
essential disks in $M(\delta)$. 
Then $Q= \hat Q\cap M$ is a punctured disk, with one boundary component, 
$\partial\hat Q$, on $\partial_1M$, and $q$ boundary components on 
$\partial_0 M$, each with slope $\delta$. 

By standard arguments we may assume that: 
\roster 
\item"(i)" $Q$ meets $T$ transversely, in properly embedded arcs and circles; 
\item"(ii)" each component of $\partial Q\cap \partial_0M$ meets each 
component of $\partial T$ in $\Delta = \Delta (\gamma,\delta)$ points; 
\item"(iii)" no arc component of $Q\cap T$ is boundary parallel in either 
$Q$ or $T$; 
\item"(iv)" no circle component of $Q\cap T$ bounds a disk in either $Q$ or $T$.
\endroster

Then, as described in \cite{GLu3, Section 2}, the arc components of 
$Q\cap T$ define graphs $G_Q$ in $\hat Q$ and $G_T$ in $\hat T$. 

The definitions and terminology of \cite{GLu1, Section 2} and 
\cite{GLu3, Section 2}, in particular the notion of a {\it $\bq$-type\/}, 
carry over to our present context. 
We assume familiarity with this terminology. 

Note that since $\hat Q$ does not separate $M(\delta)$, the signs of the 
labels around the vertices of $G_T$ (equivalently, the signs of the 
intersections of $K_\delta$ with $\hat Q$) do not necessarily alternate. 
Thus we may designate the corners of faces of $G_T$ as either $++$, $- -$, 
or $+-$, according to the signs of the corresponding pair of labels. 
However, we have the following lemma. 

\proclaim{Lemma 2.1} 
Let $\D$ be a set of disk faces of $G_T$ representing all $\bq$-types. 
Then there exists $\D' \subset \D$ such that $\D'$ represents all 
$\bq$-types and all the corners appearing in faces in $\D'$ are $+-$.
\endproclaim 

\demo{Proof} 
Let $\C_0$ be the set of $++$ and $--$ corners, and $\C_1$ the set of $+-$ 
corners, that appear in faces belonging to $\D$. 
Let $\tau_0$ be the $\C_0$-type defined by $\tau_0|(++\text{ corner})=+$, 
\break $\tau_0| (--\text{ corner}) = -$. 

Now let $\tau_1$ be any $\C_1$-type, and let $\tau$ be the 
$(\C_0\cup \C_1)$-type $(\tau_0,\tau_1)$. 
By hypothesis, there exists a face $D\in \D$ which represents $\tau$. 
We may assume, by definition of $\tau_0$, 
that the character (see \cite{GLu1, p.386}) 
of each edge endpoint at a corner of $D$ belonging to $\C_0$ is $+$. 
Since edges join points of opposite character, and since the edge endpoints 
at a corner in $\C_1$ have distinct characters, it follows that no corner 
of $D$ can belong to $\C_0$. 
We thus obtain our desired subcollection $\D'\subset \D$ such that $\D'$ 
represents all $\C_1$-types, and hence all $\bq$-types.\qed
\enddemo 

\proclaim{Theorem 2.2} 
$G_T$ does not represent all $\bq$-types.
\endproclaim 

\demo{Proof} 
Let $\D$ be a set of disk faces of $G_T$ representing all $\bq$-types. 
By Lemma~2.1, we can assume that all corners appearing in faces in $\D$ are 
$+-$, i.e., correspond to points of intersection of $K_\delta$ with $\hat Q$ 
of opposite sign. 
By \cite{GLu1, Lemma~3.1}, there exists $\D_0\subset \D$ representing 
all $\bq$-types such that each face in $\D_0$ is locally on the same side of 
$\hat Q$. 
By \cite{GLu2, Lemma~4.4}, there exists $\D_1\subset \D_0$ such that 
$\{[\partial D] : D\in \D_1\}$ is a basis for $\real^{c(\D_1)}$. 
Here $c(\D_1)$ is the number of corners appearing in faces belonging to 
$\D_1$, and $[\partial D]$ is the element of $\zed^{c(\D_1)}\subset 
\real^{c(\D_1)}$ (defined up to sign) obtained in the obvious way, by 
taking the algebraic sum of the corners in $D$. 
Tubing $\hat Q$ along the annuli in $\partial V_\delta$ corresponding 
to the corners appearing in faces belonging to $\D_1$ and surgering by the 
disks in $\D_1$ then gives a disk $\hat Q'\subset M(\delta)$ with 
$\partial \hat Q'=\partial\hat Q$, and $|\hat Q'\cap K_\delta| < 
|\hat Q\cap K_\delta| = q$, contradicting the minimality of $q$.\qed
\enddemo 

The definition of a {\it web\/} $\Lambda$ in $G_Q$ is exactly as in 
\cite{GLu3, p.601}. 
If $U$ is the component of $\hat Q-\nhd (\Lambda)$ containing $\partial \hat Q$ 
then we say $D_\Lambda = \hat Q-U$ is the {\it disk bounded by\/} $\Lambda$. 
A {\it great web\/} is a web $\Lambda$ such that $\Lambda$ contains all the 
edges of $G_Q$ that lie in $D_\Lambda$. 

Exactly as in \cite{GLu3, proof of Theorem  2.5}, Theorem~2.2 implies 

\proclaim{Theorem 2.3} 
$G_Q$ contains a great web.
\endproclaim 

Since $\Delta=3$, the arguments in \cite{GLu3, Sections 6 and 7} apply here 
without change to show

\proclaim{Theorem 2.4} 
$M(\gamma)$ does not contain a Klein bottle.
\endproclaim 

Since $M(\delta)$ is boundary-reducible and $\gamma\ne\delta$, $M(\gamma)$ 
is irreducible by \cite{S}. 

Finally we note that since $M$ is a $\que$-homology cobordism between 
$\partial_0M$ and $\partial_1M$, $M(\gamma)$ is a $\que$-homology 
$S^1\times D^2$.  
Hence $\hat T$ separates $M(\gamma)$, into $X$ and $X'$, say. 
It follows that the faces of $G_Q$ may be shaded alternately black and white, 
with the black faces lying in $X$ and the white faces lying in $X'$. 

\head 3. Scharlemann cycles and extended Scharlemann cycles\endhead 

In this section we follow \cite{GLu3, Section 3} and show that the 
relevant statements there hold in our present setting. 
For ease of reference, here and in Section~4 we will give lemmas etc.\  
the same numbers as the corresponding statements in \cite{GLu3}. 

\proclaim{Lemma 3.1} 
The edges of a Scharlemann cycle in $G_Q$ cannot lie in a disk in $\hat T$. 
\endproclaim 

\demo{Proof} 
Otherwise, as in \cite{GLu3, proof of Lemma 3.1}, $M(\gamma)$ would have a 
lens space summand, and hence be reducible.\qed
\enddemo 

\proclaim{Theorem 3.2} 
If $t\ge 4$ then $G_Q$ does not contain an extended Scharlemann cycle. 
\endproclaim 

\demo{Proof} 
The proof of \cite{GLu3, Claim 3.3} goes through verbatim in our 
present context. 

The proof of \cite{GLu3, Claim 3.4} remains valid here once we note in addition 
that $\partial M_1$ cannot be parallel to $\partial M(\gamma)$ 
($=\partial_1M$). 
For then $\hat T$ would lie in a collar of $\partial M(\gamma)$ in $M(\gamma)$, 
contradicting the essentiality of $\hat T$.

We now follow the remainder of the proof of \cite{GLu3, Theorem~3.2} 
as closely as possible. 

Let $M_2 = \overline{X-M_1}$. 
Then $X= M_1 \cup_B M_2$. 
Similarly, since $M(\gamma)$ is a $\que$-homology $S^1\times D^2$, $A$ 
separates $X'$, into $M_3$ and $M_4$, say. 

Let $T_2,T_3,T_4$ be the tori defined as $T_2,\hat T_3,\hat T_4$ in 
\cite{GLu3, pp.609--610}. 
These satisfy $|T_i\cap K_\gamma| < |\hat T\cap K_\gamma|$, $i=2,3,4$, and 
hence, by the minimality of $t$, $T_i$ is either compressible, and therefore 
(since $M(\gamma)$ is irreducible) bounds a solid torus in $M(\gamma)$, or 
is peripheral, i.e., parallel to $\partial M(\gamma)$, $i=2,3,4$.  
It follows that exactly one of $M_2,M_3,M_4$ is a collar of 
$\partial M(\gamma)$ in $M(\gamma)$, and the other two are solid tori. 

Now let $T_0$ be the torus $A\cup (C_1-A_1)\cup B$. 
Then $T_0$ separates $M(\gamma)$ into $Y= M_1\cup_A M_3$ and 
$Y'= M_2\cup_{C_2} M_4$. 
Also $|T_0\cap K_\gamma| <t$. 
Hence again $T_0$ either bounds a solid torus or is peripheral. 

If $T_0$ is peripheral, then $M(\gamma) \cong Y$ or $Y'$, the union of two 
solid tori along an annulus essential in the boundary of each. 
Since $M(\gamma)$ is irreducible, this implies that $M(\gamma)$ is a Seifert 
fiber space over the disk with at most two singular fibers, contradicting 
the fact that it contains an essential torus. 

If $Y$, say, is a solid torus, then $M_2$ or $M_4$ (say $M_2$) is a collar 
of $\partial M(\gamma)$, and hence $M(\gamma) \cong Y\cup M_4$ is a 
union of two solid tori, which is a contradiction as before. 
We get a similar contradiction if $Y'$ is a solid torus.\qed
\enddemo 

We remark that we will only use Theorem 3.2 for extended Scharlemann cycles 
of length~2; see Section~4. 
However, we have given the proof of the general statement, for possible 
future reference. 

Since $M(\gamma)$ does not contain a Klein bottle by Theorem~2.4, the proof 
of \cite{GLu3, Lemma~3.10} applies verbatim to give the same statement here: 

\proclaim{Lemma 3.10} 
Let $\sigma_1$ and $\sigma_2$ be Scharlemann cycles of length $2$ in $G_Q$ 
on distinct label-pairs.  
Then the loops on $\hat T$ formed by the edges of $\sigma_1$ and $\sigma_2$ 
respectively are not isotopic on $\hat T$. 
\endproclaim 

Similarly, the proof of \cite{GLu3, Theorem 3.11} gives 

\proclaim{Theorem 3.11} 
At most three labels can be labels of Scharlemann cycles of length $2$ in 
$G_Q$. 
\endproclaim 

\head 4. The case $t\ge 4$\endhead 

We obtain a contradiction in this case exactly as in \cite{GLu3, Section~4}. 

Let $\Lambda$ be a great web in $G_Q$, as guaranteed by Theorem~2.3, and 
$D_\Lambda$ be the disk that it bounds. 
For every label $x$ of $G_Q$ let $\Lambda_x$ be the graph in $D_\Lambda$ 
consisting of the vertices of $\Lambda$ along with all $x$-edges in $\Lambda$. 

Let $V$ be the number of vertices of $\Lambda$. 

\proclaim{Lemma 4.2} 
Let $x$ be a label of $G_Q$. 
If $\Lambda_x$ has at least $3V-4$ edges, then $\Lambda_x$ contains a bigon.
\endproclaim 

\demo{Proof} 
The proof of \cite{GLu3, Lemma 4.2} applies here without change. 

We take this opportunity to thank Chuichiro Hayashi and Kimihiko Motegi 
for pointing out that the case where the face $f$ is a monogon was not 
explicitly addressed in \cite{GLu3}. 
However, in that case $\Lambda_x$ could have at most one ghost $x$-label, 
since otherwise $\Lambda$ would have more than $t$ ghost labels. 
Hence $E\ge 3V-1$. 
Also, if $\Lambda_x$ does not contain a bigon then $1+3(F-1) \le 2E$. 
Therefore $2=V-E+F\le (E+1) /3 - E+ 2(E+1)/3=1$, a contradiction.\qed
\enddemo 

Lemma 4.2 and Theorems 3.2 and 3.11 then lead to a contradiction exactly 
as in \cite{GLu3, p.615}. 

\head 5. The case $t=2$\endhead 

As in the previous section, let $\Lambda$ be a great web in $G_Q$. 

\proclaim{Lemma 5.1} 
$\Lambda$ contains a bigon.
\endproclaim 

\demo{Proof} 
Regard $\Lambda$ as a graph in the disk, with $V$ vertices, $E$ edges, and 
$F$ faces. 

Since $\Lambda$ has at most two ghost labels, $2E\ge 6V-2$, so $E\ge3V-1$. 
Suppose $\Lambda$ contains no bigon; then $3F < 2E$. 
Therefore $1= V-E+F< (E+1)/3 - E+ 2E/3 = \frac13$, a contradiction.\qed
\enddemo 

By the parity rule, each edge of $\Lambda$ must connect different vertices 
of $G_T$ when considered as edges lying in $G_T$. 
Recall that there are four {\it edge classes\/} in $G_T$, i.e., isotopy 
classes of non-loop 
edges of $G_T$ in $\hat T$ rel$\{$vertices of $G_T\}$, which we 
call $1,\alpha,\beta,\alpha\beta$, as illustrated in \cite{GLu3, Figure~7.1}. 
The ordering of these classes around vertices~1 and 2 of $G_T$ is indicated 
in Figure~5.1. 
We may label an edge $e$ of $G_Q$ by the class of the corresponding edge 
of $G_T$; we refer to this label as the {\it edge class label\/} of $e$. 

\fig{1.5}{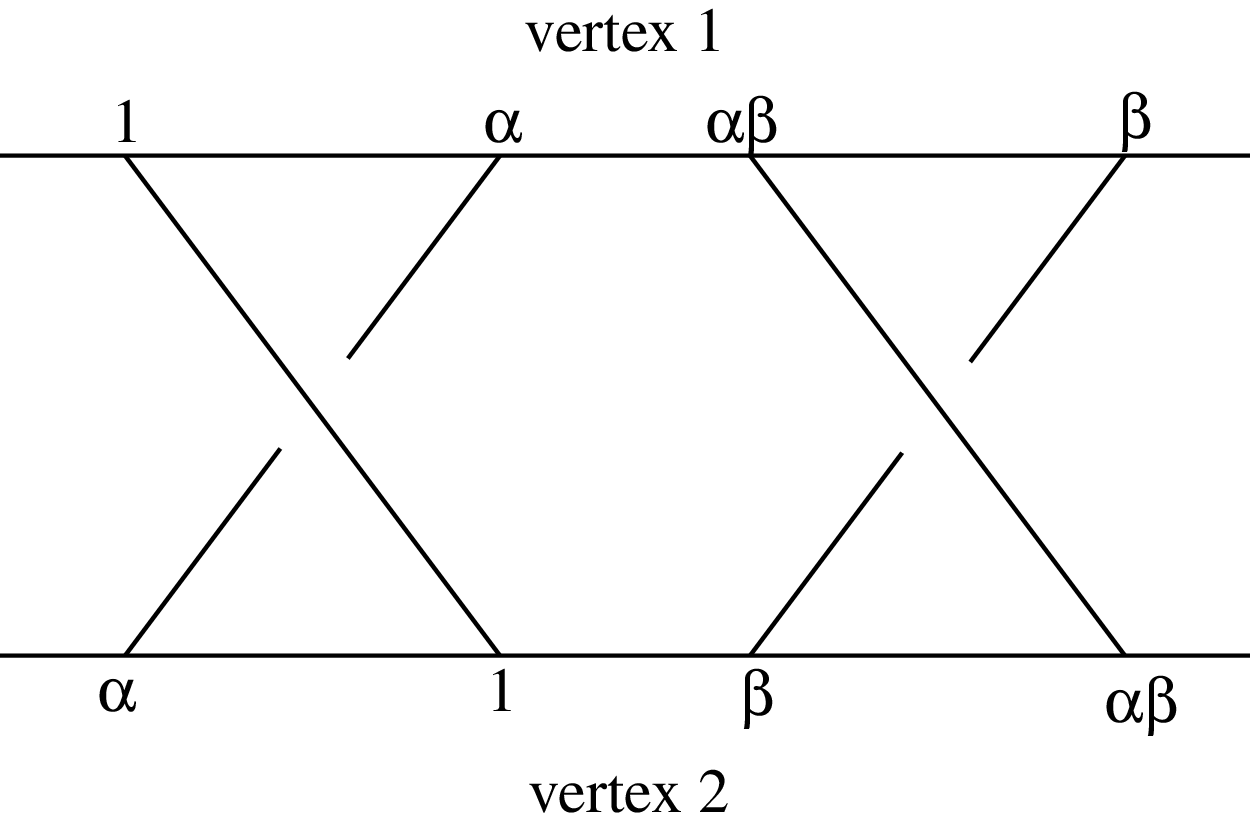}{5.1} 

\proclaim{Lemma 5.2}
Any two black (white) bigons in $\Lambda$ have the same pair of edge class 
labels. 
\endproclaim 

\demo{Proof} 
Let $f,f'$ be two (say) black bigons in $\Lambda$, with vertices $x,x'$ 
and $y,y'$, and edge class labels $\lambda,\mu$ and $\lambda',\mu'$, 
respectively. 
See Figure~5.2. 
Note that $\lambda\ne\mu$ and $\lambda'\ne \mu'$ by Lemma~3.1. 

\fig{.75}{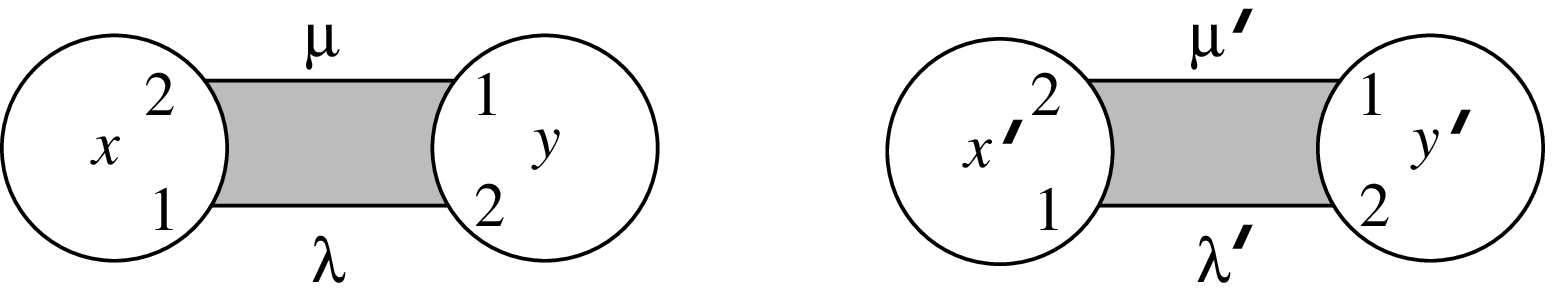}{5.2} 

Let $H_{12}$ be the 1-handle $V_\gamma\cap X$. 
The anticlockwise ordering of the labels $x,x',y,y'$ around vertex~1 of $G_T$ 
must agree with their clockwise ordering around vertex~2, since the corners 
of $f$ and $f'$ give four disjoint arcs on the annulus 
$\overline{\partial H_{12}-\hat T}$ joining corresponding labels. 
However, if $\lambda=\lambda'$ but $\mu\ne \mu'$ then it is easy to see 
that these orderings do not agree. 
For example, the case $\lambda=\lambda'=1$, $\mu=\alpha$, $\mu'=\beta$ 
is illustrated in Figure~5.3. 
Similarly, the pairs $\{1,\alpha\beta\}$, $\{\alpha,\beta\}$ give 
incompatible orderings of the labels. 
\smallskip

\fig{1.5}{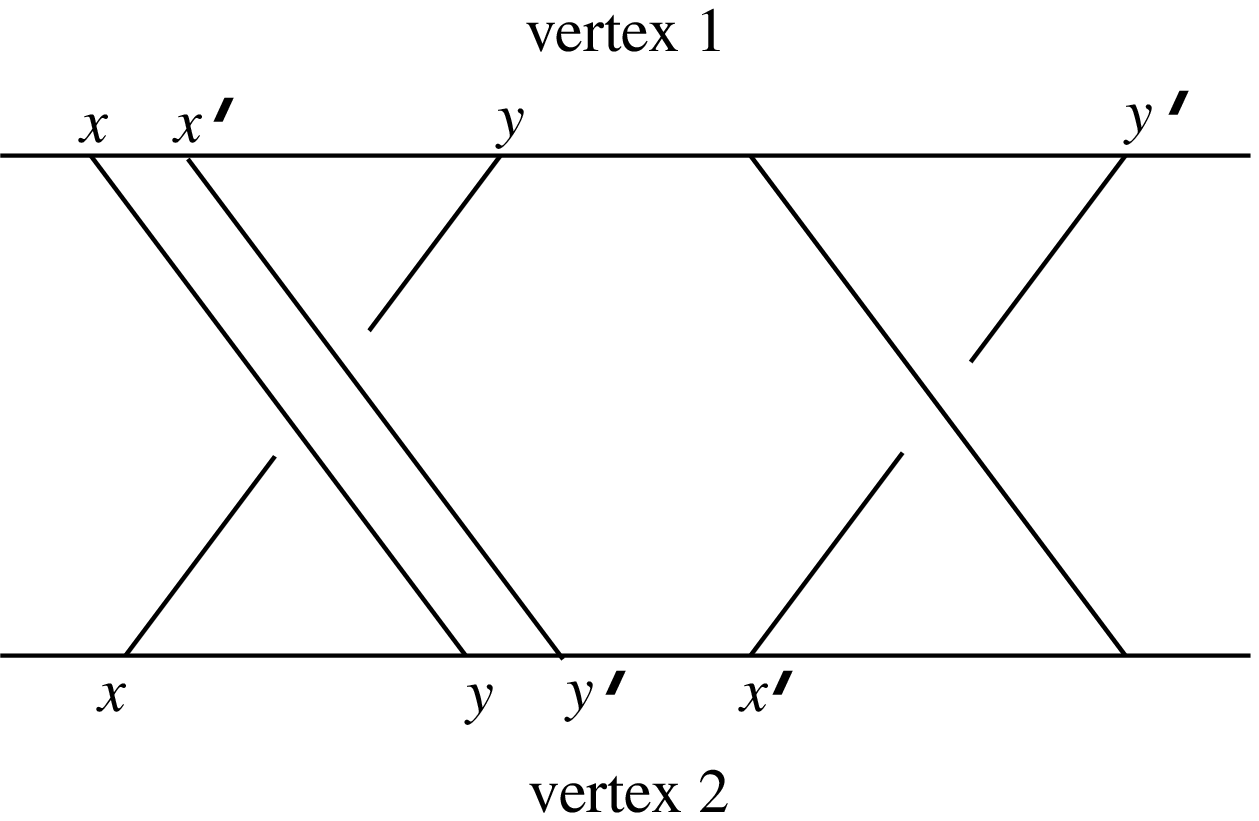}{5.3} 

Thus the only possible distinct pairs of edge class labels for $f$ and $f'$  
are $\{1,\alpha\}$ and $\{\alpha\beta,\beta\}$ (or $\{1,\beta\}$ and 
$\{\alpha\beta,\alpha\}$). 
Now by shrinking $H_{12}$ to its core, $H_{12}\cup f$ becomes a M\"obius 
band $B$ in $X$ such that $\partial B$ is the loop on $\hat T$  formed by 
the edges of $f$. 
Similarly $f'$ gives rise to a M\"obius band $B'$. 
If the edge class labels for $f$ and $f'$ are as above, then $\partial B$ 
and $\partial B'$ may be isotoped on $\hat T$ to be disjoint. 
This produces a Klein bottle in $M(\gamma)$, contradicting Theorem~2.4.\qed
\enddemo 

\proclaim{Lemma 5.3} 
Let $e,e'$ be edges of $\Lambda$ incident to a vertex $v$ of $\Lambda$, with 
the same label at $v$. 
Then $e$ and $e'$ have distinct edge class labels. 
\endproclaim 

\demo{Proof} 
If not, then $e$ and $e'$ would be parallel in $G_T$, and hence would cobound 
a family of $q+1$ parallel edges of $G_T$. 
The argument of \cite{GLi, p.130, Case~(2)} now constructs a cable space in 
$M$, contradicting the assumption that $M$ is simple.\qed
\enddemo 

If $f$ is a white (black) face of $\Lambda$ such that each edge of $f$ 
belongs to a black (resp. white) bigon, then we say that $f$ is 
{\it surrounded by bigons\/}. 

\proclaim{Lemma 5.4} 
$\Lambda$ does not contain a face surrounded by bigons.
\endproclaim 

\demo{Proof} 
Suppose that $\Lambda$ contains, say, a white face $f$ surrounded by black 
bigons. 

By Lemma 5.2 (and Lemma 3.1) all the black bigons have the same pair of 
(distinct) edge class labels $\lambda,\mu$. 

Let $v$ be a vertex of $f$. 
By Lemma~5.3 the edges of the two black bigons incident to $v$ having 
label~1 (say) at $v$ have distinct edge class labels. 
Therefore the two edges of $f$ incident to $v$ have the same edge class label 
(see Figure~5.4). 
Hence all the edges of $f$ have the same edge class label, contradicting 
Lemma~3.1.\qed
\smallskip
\fig{1.0}{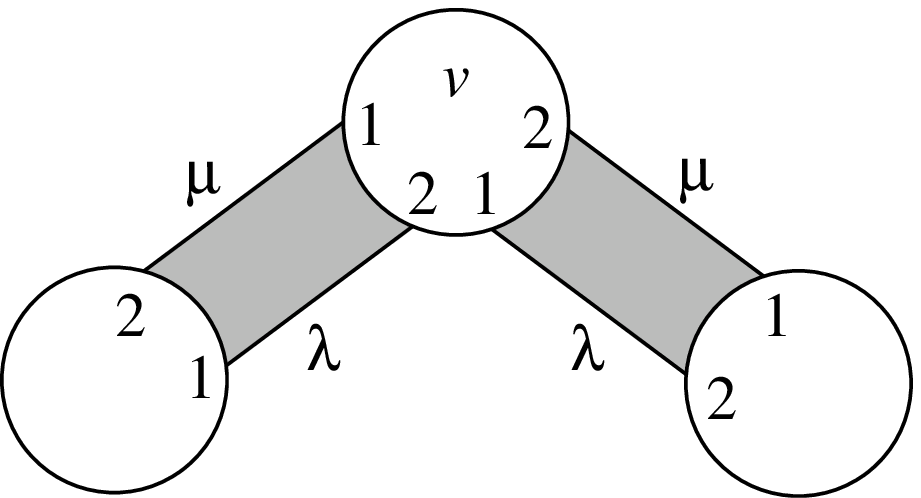}{5.4} 
\enddemo 

\proclaim{Lemma 5.5} 
Not all black faces of $\Lambda$ are bigons. 
Similarly, not all white faces of $\Lambda$ are bigons.
\endproclaim 

\demo{Proof} 
Suppose, for example, that all black faces of $\Lambda$ are bigons. 

Define a graph $\Sigma$ in the disk $\hat Q$ as follows. 
The vertices of $\Sigma$ are the vertices of $\Lambda$, and the edges 
of $\Sigma$ are in one-one correspondence with the black bigons of 
$\Lambda$, each edge joining the vertices at the corners of the 
corresponding bigon and lying in its interior. 
Then each vertex of $\Sigma$ has valency at least~2, except that if 
$\Lambda$ has a vertex at which there are two ghost labels, then 
that vertex may have valency~1 in $\Sigma$. 
Hence if $V$ and $E$ are the number of vertices and edges of $\Sigma$ 
respectively, then $2E \ge 2V-1$, and so $2E\ge 2V$. 
Therefore 
$$\align 
1 &\le V-E+\sum \chix (f) \text{ (the sum being over all faces $f$ 
of $\Sigma$)} \cr
& \le \sum\chix (f)\ ,
\endalign$$
implying that $\Sigma$ has a disk face. 
But such a face corresponds to a white face of $\Lambda$ surrounded by 
black bigons, contradicting Lemma~5.4.\qed
\enddemo

We shall say that two faces $f_1,f_2$ of $\Lambda$ of the same color are 
{\it isomorphic\/} if the cyclic sequences of edge class labels obtained 
by reading around the boundaries of $f_1$ and $f_2$ in the same direction,  
are equal. 

\proclaim{Lemma 5.6} 
Either all black faces of $\Lambda$ are isomorphic or all white faces 
of $\Lambda$ are isomorphic. 
\endproclaim 

\demo{Proof} 
Recall that $\hat T$ separates $M(\gamma)$, into $X$ and $X'$. 
We suppose without loss of generality that $\partial M(\gamma) \subset X$. 

Recall also that 
$H_{12}$ is the 1-handle $V_\gamma\cap X$. 
Let $f$ be any black face of $\Lambda$, and let $N= \nhd (\hat T\cup H_{12} 
\cup f) \subset X$. 
Then (since all the vertices of $\Lambda$ have the same sign) $\partial N 
= \hat T\cup T_1$, where $T_1$ is a torus. 
Since $T_1 \cap K_\gamma = \emptyset$, and since $M$ is irreducible and 
atoroidal, $T_1$ is parallel to $\partial_1 M = \partial M(\gamma)$. 
Hence $W= \overline{X-H_{12}}$ is a compression body, with $\partial_+W$ 
the genus~2 surface obtained by adding the tube $\overline{\partial H_{12} 
- \hat T}$ to $T$, and $\partial_-W = \partial M(\gamma)$. 
It follows that  $f$ is the unique non-separating disk in $W$,   
up to isotopy. 

Let $f'$ be any other black face of $\Lambda$. 
Then $\partial f$ and $\partial f'$ are isotopic in $\partial_+ W$, 
and hence (freely) homotopic in $\hat T\cup H_{12}$. 
Now $\pi_1(\hat T\cup H_{12}) \cong \pi_1(\hat T) * \zed$, where, taking 
as base-``point'' a disk neighborhood in $\hat T$ of an edge in $G_T$ 
in edge class~1, $\pi_1 (\hat T) \cong \zed\times\zed$ has basis 
$\{\alpha,\beta\}$, represented by edges in the correspondingly named 
edge classes, oriented from vertex~2 to vertex~1, and $\zed$ is generated 
by $x$, say, represented by an arc in $H_{12}$ going from vertex~1 to 
vertex~2. 
Then, if the sequence of edge class labels around $\partial f$ (in the 
appropriate direction) is $(\gamma_1,\ldots,\gamma_n)$, $\partial f$ 
represents $\gamma_1 x\gamma_2 x\ldots \gamma_n x\in \pi_1(\hat T) *\zed$, 
and similarly for $\partial f'$. 
Since $\partial f$ and $\partial f'$ are homotopic in $\hat T\cup H_{12}$, 
we conclude that the corresponding cyclic sequences $(\gamma_1,\ldots, 
\gamma_n)$ and $(\gamma'_1,\ldots,\gamma'_{n'})$ are equal, i.e., $f$ and 
$f'$ are isomorphic. 
Thus all black faces of $\Lambda$ are isomorphic, and the lemma is proved.\qed 
\enddemo 

By Lemma 5.1, $\Lambda$ contains a bigon. 
We may therefore assume from now on that $\Lambda$ contains a black bigon. 

\proclaim{Lemma 5.7} 
All white faces of $\Lambda $ are isomorphic.
\endproclaim 

\demo{Proof} 
If not, then by Lemma 5.6 all black faces of $\Lambda$ are isomorphic, 
and therefore bigons. 
But this contradicts Lemma~5.5.\qed
\enddemo 

\proclaim{Lemma 5.8} 
$\Lambda$ does not contain a white bigon.
\endproclaim 

\demo{Proof} 
This follows from Lemmas 5.7 and 5.5.\qed
\enddemo 

Recall that $D_\Lambda$ is the disk bounded by the great web $\Lambda$. 
Travelling around $\partial D_\Lambda$ in some direction defines in the 
obvious way a cyclic sequence of edges of $\Lambda$, $(e_1,\ldots ,e_n)$. 
(Note that the same edge may appear twice, with opposite orientations.) 
If the subgraph of $\Lambda$ consisting of these edges has no cut vertex, 
define $\Lambda_0=\Lambda$. 
(Recall that a cut vertex of a connected graph is one whose removal 
disconnects the graph.) 
Otherwise, let $v_0$ be an outermost cut vertex of this subgraph. 
This vertex cuts off a subdisk of $D_\Lambda$ which is the 
disk $D_{\Lambda_0}$ bounded 
by a subgraph $\Lambda_0$ of $\Lambda$ (see Figure~5.5). 
We shall call $\Lambda_0$ an {\it extremal subgraph\/} of $\Lambda$. 
If $\Lambda_0\ne \Lambda$ we call $v_0$ the {\it attaching vertex\/} 
of $\Lambda_0$. 
\smallskip

\fig{1.75}{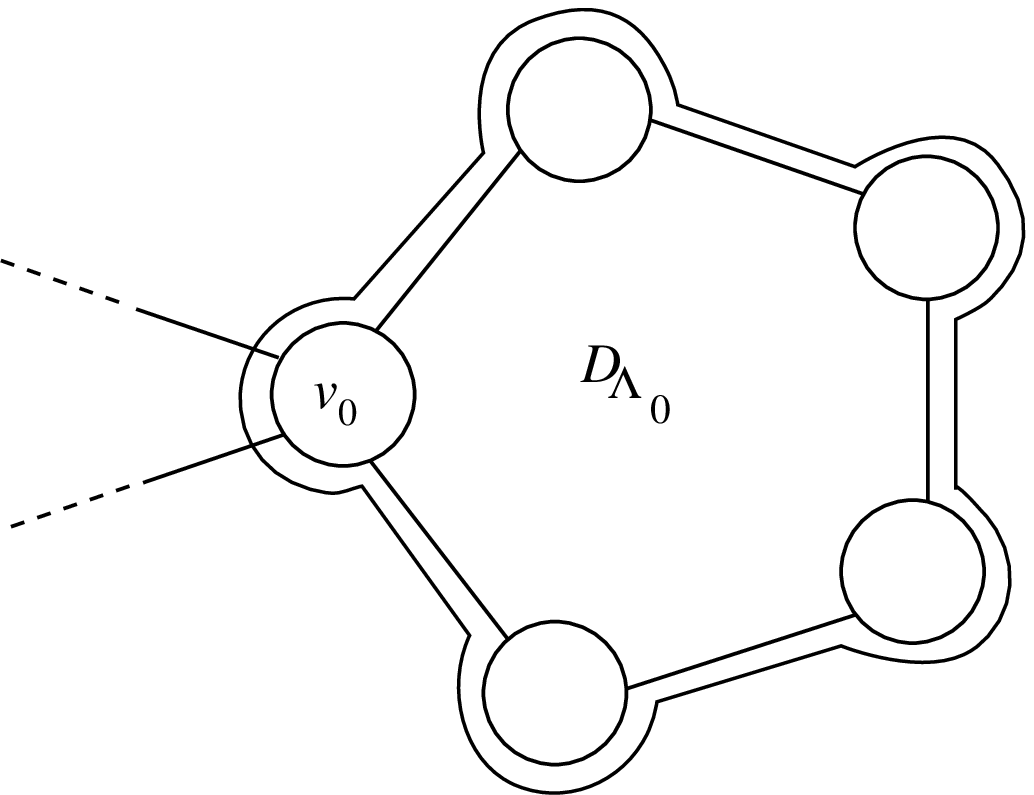}{5.5} 

An {\it interior vertex\/} of $\Lambda$ is a vertex of $\Lambda$ that 
is not an endpoint of any of the edges $e_1,\ldots,e_n$. 
All six corners at such a vertex belong to faces of $\Lambda$. 

\proclaim{Lemma 5.9} 
$\Lambda$ contains an interior vertex. 
\endproclaim 

\demo{Proof} 
The following proof was suggested by the referee. 

Since $\Lambda$ has at most two ghost labels, $\Lambda$ has an extremal 
subgraph $\Lambda_0$ such that either $\Lambda_0 = \Lambda$ or $\Lambda_0$ 
contains at most one ghost label of $\Lambda$ at a vertex other than the 
attaching vertex $v_0$ of $\Lambda_0$. 
If $\Lambda=\Lambda_0$, we can also choose a vertex $v_0$ of $\Lambda$ 
such that $\Lambda$ has at most one ghost label at a vertex other than $v_0$. 
We shall show that $\Lambda_0$ contains an interior vertex. 

Suppose for contradiction that $\Lambda_0$ has no interior vertices. 
Let $m$ be the number of vertices of $\Lambda_0$. 
Note that $m>1$, since $G_Q$ has no trivial loops. 
Let $\oLambda$ be the reduced graph of $\Lambda_0$, obtained 
by amalgamating each family of parallel edges to a single edge. 
If $m=2$ or 3, then each vertex of $\oLambda$ has valency~1 or 2 
respectively. 
If $m\ge4$, the interior edges of $\oLambda$ decompose 
$D_{\Lambda_0}$ into subdisks, and by considering outermost such subdisks 
we see that $\oLambda$ has at least two vertices of valency~2. 
Hence in all cases $\oLambda$ has a vertex $v\ne v_0$ of valency~1 or 2. 
Since there is at most one ghost label at $v$, $\Lambda_0$ has a family 
of at least three parallel edges at $v$, and hence contains a white bigon, 
contradicting Lemma~5.8.\qed 
\enddemo 

\proclaim{Lemma 5.10} 
Each white face of $\Lambda$ has at least three distinct edge class labels. 
\endproclaim 

\demo{Proof} 
By Lemma 5.9, $\Lambda$ contains an interior vertex $v$. 
Then all three white corners at $v$ belong to faces of $\Lambda$. 
By Lemma~5.3, the three edges with label (say) 1 at $v$ have distinct 
edge class labels. 
The result now follows from Lemma~5.7.\qed
\enddemo 

\proclaim{Lemma 5.11} 
Each white face of $\Lambda$ has length at least $4$. 
\endproclaim 

\demo{Proof} 
The argument in the first paragraph of \cite{GLu3, proof of Lemma~3.7} 
shows that a face of $\Lambda$ of length~3 has only two distinct edge 
class labels, contradicting Lemma~5.10.\qed
\enddemo 

Recall that $\Lambda$ has either 0 or 2 ghost labels. 
In the former case, let $\Lambda^+=\Lambda$. 
In the latter case, let $\Lambda^+\subset \hat Q$ be obtained from 
$\Lambda$ by adjoining two disjoint arcs running from the ghost labels 
to $\partial \hat Q$. 
The complementary regions of $\Lambda^+$ in $\hat Q$ are the faces of 
$\Lambda$ together with either one or two {\it outside regions\/}, i.e., 
regions that meet $\partial \hat Q$. 
The black/white shading of the faces of $\Lambda$ extends to a shading 
of the complementary regions of $\Lambda^+$. 

Define a graph $\Gamma$ in the disk $\hat Q$ as follows. 
The vertices of $\Gamma$ are the ({\it fat\/}) vertices of $\Lambda$, 
together with {\it dual vertices\/}, one in the interior of each black 
complementary region of $\Lambda^+$ that is either a face of $\Lambda$ 
of length at least~3, or an outside region. 
The edges of $\Gamma$ consist of edges joining each dual vertex to the 
fat vertices in the boundary of the corresponding region, together with 
an edge in the interior of each black bigon in $\Lambda$, joining the 
(fat) vertices at the ends of the bigon. 
(If $\Lambda^+ =\Lambda$ and the outside region is black, then there is 
a choice involved in defining the edges incident to the dual vertex 
corresponding to that region.) 

\proclaim{Lemma 5.12} 
$\Gamma$ has a face of length at most $5$.
\endproclaim 

\demo{Proof} 
The valency of each fat vertex of $\Gamma$ is 3. 
The valency of each dual vertex corresponding to a face of $\Lambda$ 
is at least~3. 
There is at most one dual vertex corresponding to an outside region, and 
that vertex has valency at least~1. 
Hence, if $\Gamma$ has $V$ vertices and $E$ edges then $2E\ge 3V-2$. 
If $\Gamma$ has $F$ faces, and each face has length at least~6, then $6F< 2E$. 
Therefore $1= V-E+F < 2(E+1)/3 -E+E/3 = \frac23$, a contradiction.\qed
\enddemo 

Let $g$ be a face of $\Gamma$ as in Lemma 5.12. 
Then $g$ contains a unique white face $f$ of $\Lambda$. 
If $f$ has $b$ edges which are edges belonging to black bigons of $\Lambda$, 
and $c$ other edges, then the length of $g$ is $b+2c\le 5$. 
On the other hand $b+c \ge4$ by Lemma~5.11. 
Also, $c>0$ by Lemma~5.4; hence $c=1$. 
Thus every edge of $f$ except one belongs to a black bigon, and the 
argument in the proof of Lemma~5.4 shows that every such edge has the 
same edge class label. 
Therefore $f$ has at most two edge class labels, contradicting Lemma~5.10.  

We have thus shown that the case $t=2$ is impossible. 

This completes the proof of Theorem 1.1. 

\head 6. Knots in solid tori\endhead 

In \cite{GLu3} and \cite{GLu4} it is shown that if $M$ is simple, 
$M(\gamma)$ is toroidal, $M(\delta)$ is $S^3$, and $\Delta (\gamma,\delta) 
= 2$, then $t=2$. 
\cite{MM} gives examples of  simple manifolds, $M$, for which $M(\gamma)$ 
is toroidal, $M(\delta) = S^1\times D^2$, and $\Delta (\gamma,\delta)=2$. 
In these examples $t=2$. 
It is natural then to ask if this is always the case. 

\proclaim{Theorem 6.1} 
Let $M$ be simple 3-manifold such that $M(\gamma)$ is toroidal, 
$M(\delta)$ is $S^1\times D^2$, and $\Delta (\gamma,\delta)=2$.  
Let $K_\gamma$ be the core of  the attached solid torus in $M(\gamma)$ 
and $\hat T$ be an essential torus in $M(\gamma)$ 
that intersects $K_\gamma$ minimally. 
Then $|K_\gamma \cap \hat T| =2$. 
\endproclaim 

The remainder of this section is devoted to the proof of Theorem~6.1. 
Let $\lambda =S^1\times *$, 
$\mu = *\times \partial D^2$ in $\partial_1M = \partial M(\delta)$ 
be a longitude and meridian respectively of $S^1\times D^2$. 

For any integer $n$, let $M_n$ be the $(\lambda +n\mu)$-filling of $M$. 

\proclaim{Lemma 6.2}
$M_n$ is simple for $n$ sufficiently large.
\endproclaim 

\demo{Proof} 
By \cite{Th}, $M$ is hyperbolic and all but finitely many fillings on a 
single component of $\partial M$ are hyperbolic. 
See also Theorem~1.3 of \cite{G}.\qed
\enddemo 

Let $N=M(\gamma)$, and parametrize slopes (as rational numbers) 
on $\partial N= \partial_1M$ using $\lambda,\mu$. 
Thus $M_n(\gamma) = N(n)$. 

Let $\hat T$ be an essential torus in $N$ that minimizes 
$|K_\gamma\cap \hat T|$. 

\proclaim{Lemma 6.3} 
If $\hat T$ compresses in $N(n_1)$ and $N(n_2)$ then either 
$|K_\gamma\cap \hat T|=2$ or $|n_1-n_2|\le 2$. 
\endproclaim 

\demo{Proof} 
We assume for contradiction that every essential torus in $N$ (which must 
separate since $\Delta (\gamma,\delta)=2$) intersects $K_\gamma$ at least 
four times and that $\hat T$ compresses in $N(n_1),N(n_2)$ 
where $|n_1-n_2|\ge3$. 
Applying \cite{CGLS, Theorem 2.4.2} to the side of $\hat T$ containing 
$\partial_1 M$, 
this last condition guarantees an essential annulus $A'\subset N$ such 
that one boundary component of $A'$ is $A'\cap \hat T$, the other is 
$A'\cap \partial_1 M$, and  
$A'\cap \partial_1 M$ is distance~1 from each of the integral slopes $n_1$ 
and $n_2$. 
Since $|n_1-n_2|\ge3$, $A'\cap \partial_1 M$ must be a copy of $\mu$. 
Surgering $\hat T$ along $A'$ gives a properly embedded annulus $A''\subset N$ 
whose boundary consists of two copies of $\mu$. 
Since $\hat T$ is essential in $N$, so is $A''$. 

Let $A$ be an essential annulus in $N$, with boundary consisting of 
two copies of $\mu$, such that $|K_\gamma \cap A|$ is minimal over all 
such annuli.  
Since $M$ contains no essential annulus, 
$K_\gamma \cap A\ne \emptyset$. 
Let $D$ be a  meridian disk of $M(\delta)= S^1\times D^2$, chosen so as to 
minimize $|K_\delta \cap D|$. 
We may isotope $D$ so that $\partial D\cap \partial A=\emptyset$. 
Then $T'=A\cap M$ and $Q=D\cap M$ are essential planar surfaces in $M$, 
whose ``outer'' boundaries, $\partial A$ and 
$\partial D$ lie on $\partial_1 M$ and are disjoint, and whose ``inner'' 
boundaries $\partial T' \cap \partial_0 M$, $\partial Q\cap 
\partial_0 M$ consist of, say, $t'$ 
copies of $\gamma$ and $q$ copies of $\delta$ respectively. 

Recall that $\Delta (\gamma,\delta)=2$. 
As in Section~2 the arc components of $T'\cap Q$ give rise to graphs 
$G_{T'},G_Q$ in $A,D$ (resp.). 
By abstractly identifying the two boundary components of $A$, we may regard 
$G_{T'}$ as a graph in a torus, and we are now exactly 
in the combinatorial set-up of Section~2. 
In particular we may apply Theorem~2.3 to conclude that $G_Q$ contains a 
great web, $\Lambda$. 
(Note that the one face of $G_{T'}$ in the torus which is not 
a face of $G_{T'}$ in $A$ is not a disk face, hence would not be involved 
in any collection of faces representing all types). 
Because $\Lambda$ has at most $t'$ ghost labels and $\Delta(\gamma,\delta)=2$, 
there is a label $x$ of $G_Q$ and a vertex   $y_0$ of $\Lambda$ such 
that the following conditions are satisfied: 
\roster
\item"(i)" for any vertex $y$ of $\Lambda$ other than $y_0$, there is an 
edge of $\Lambda$ incident   to $y$ at each occurrence of the label $x$ 
on $y$, and 
\item"(ii)" there is an edge of $\Lambda$ incident to $y_0$ at some 
occurrence of the label $x$ at $y_0$.
\endroster
If we let $\Sigma$ be the set of all $x$-edges of $\Lambda$, then $\Sigma$ 
is a {\it great $x$-web\/} in the sense of \cite{GLu2, p.390}. 
The argument of Theorem~2.3 of \cite{GLu2} now shows that $A$ is separating 
and either (a) $G_Q$ must contain Scharlemann cycles on distinct label pairs 
or (b) $t'=2$. 
(This argument uses Lemma~2.2 of the same paper, which requires the fact that 
consecutive labels on $G_Q$ represent vertices of $G_{T'}$ of opposite sign. 
This  assumption was shown not to be necessary in 
\cite{HM2, Proposition~5.1}. 
However, one may guarantee this condition 
in the present context by showing that $A$ 
is separating. 
This is done in the same way as in Theorem~2.3: show that $\Lambda$ contains 
a Scharlemann cycle, then use the face of $\Lambda$ bounded by this 
Scharlemann cycle to tube and compress $A$ to get a new annulus $A'$ such 
that $|A'\cap K_\delta| < |A\cap K_\delta|$. 
If $A$ were non-separating, $A'$ would be also, contradicting the minimality 
of $A$.)

Lemma 3.1 of \cite{HM1} shows that conclusion (a) contradicts the minimality 
of $|A\cap K_\gamma|$. 
Thus we must have that $t'=2$. 
Let $X_1$ and $X_2$ be the components into which $A$ separates $N$. 
Then each $\partial X_i$ is a torus meeting $K_\gamma$ twice, and hence 
is compressible. 
Since $N$ is irreducible by \cite{S}, $X_i$ must be a solid torus. 
As $N= X_1\cup_A X_2$, $N$ is Seifert-fibered over the disk with at most 
two exceptional fibers. 
But this contradicts the fact that $N$ contains an essential torus. 
This contradiction finishes the proof of Lemma~6.3.\qed
\enddemo 

\remark{Remark} 
In fact, the possibility that $N$ contains an essential torus intersecting 
$K_\gamma$ twice can be excluded from the conclusion of Lemma~6.3. 
For example one can argue as above to show the existence of $A$, 
then argue as in \cite{HS} to arrive at a contradiction. 
In particular, Sections~4 and 5 of that paper along with Lemma~3.1 of 
\cite{HM1} show that the great $x$-web $\Lambda$ that we see in the 
above argument cannot occur. 
\endremark 

\demo{Proof of Theorem 6.1} 
Lemmas 6.2 and 6.3 show that if $|K_\gamma \cap \hat T|\ne 2$, then 
there is an integer $J$ such that for all $n>J$, $M_n$ is a simple 
manifold and $M_n(\gamma)$ contains an essential torus. 
On the other hand $M_n(\delta)=S^3$, so $M_n$ is the exterior of a simple 
knot in $S^3$. 
By Theorem~1.2 of \cite{GLu3}, (which is proved in \cite{GLu3} and 
\cite{GLu4}), each $M_n$ contains an essential, 
twice-punctured torus, $T_n$, where $\partial T_n$ is two copies of $\gamma$. 
If, for some $n$, $T_n$ can be isotoped to lie inside $M$, then we are done. 
So assume not. 
Then for each $n>J$, there is an incompressible, $\partial$-incompressible 
surface $T'_n$ in $M$ whose boundary is a non-empty collection of copies 
of $\lambda +n\mu$ on $\partial_1 M$ and two copies 
of $\gamma$ on  $\partial_0 M$. 
But this contradicts \cite{H}.\qed
\enddemo 

In the above proof of Theorem 6.1, the direct appeal to Theorem~1.2 of 
\cite{GLu3} was possible because of the assumption that $M(\delta)$ is 
a solid torus. 
Conceivably Theorem~6.1 is still true if it is only assumed that $M(\delta)$ 
is boundary-reducible. 
If so, a proof in this generality might be obtainable by suitably modifying 
the proof, in \cite{GLu3} and \cite{GLu4}, of Theorem~1.2 of \cite{GLu3}.

\enddocument